\documentclass[11pt]{article}
\usepackage{amsmath}
\usepackage{amssymb}

\usepackage{amsmath}
\usepackage{amssymb}
\usepackage{latexsym}
\usepackage{epsfig}

{\catcode `\@=11 \global\let\AddToReset=\@addtoreset}
\AddToReset{equation}{section}

\usepackage{epsfig}
\usepackage{amsmath}
\usepackage{amssymb}
\usepackage{graphicx}
\usepackage{color}

\usepackage[latin1]{inputenc}

\def\I{\mathbb{I}}

\newtheorem{theorem}{Theorem}
\newtheorem{lemma}{\bf Lemma}
\newtheorem{corollary}{\bf Corollary}
\newtheorem{proposition}{Proposition}
\newtheorem{@definition}{\sc Definition}

\newtheorem{@remark}{\sc Remark}

\newtheorem{@example}{\sc Example}
\newenvironment{example}{\begin{@example}\rm}{\end{@example}}

\newtheorem{thm}{Theorem}

\newcommand{\beqn}{\begin{displaymath}}
\newcommand{\eeqn}{\end{displaymath}}
\newcommand{\beq}{\begin{equation}}  
\newcommand{\eeq}{\end{equation}}

\def\mathsf{\bf}
\def\N{\mathbb{N}}

\def\R{\mathbb{R}}
\def\C{\mathbb{C}}
\def\Z{\mathbb{Z}}

\def\I{\mathbb{I}}
\def\i{\mathrm i}
\def\d{\mathrm d}
\def\e{\mathrm e}

\def\E{\mathrm E}
\def\P{\mathrm P}

\def\text{\mbox}

\def\1{{\bf 1}}

\newcommand{\Leb}{{\rm Leb}}

\newcommand{\nn}{\nonumber}
\newcommand{\noi}{\noindent}

\newcommand{\mbt}{\boldsymbol{t}}

\newcommand{\mbu}{\boldsymbol{u}}

\newcommand{\mbx}{\boldsymbol{x}}

\newcommand{\mbz}{\boldsymbol{z}}

\newcommand{\mby}{\boldsymbol{y}}

\newcommand{\mbgamma}{\boldsymbol {\gamma}}

\definecolor{vp}{rgb}{0.55, 0.71, 0.0}

\def\eq2{
\stackrel{\small \rm mod \,2}{=}}

\def\n2{
\stackrel{\small \rm mod \,2}{\neq}}

\def\limd{\renewcommand{\arraystretch}{0.5}
\begin{array}[t]{c}
\stackrel{\rm d}{\longrightarrow} \\
\end{array}\renewcommand{\arraystretch}{1}}

\def\limfdd{\renewcommand{\arraystretch}{0.5}
\begin{array}[t]{c}
\stackrel{\rm fdd}{\longrightarrow} \\
\end{array}\renewcommand{\arraystretch}{1}}

\newtheorem{rem}[thm]{Remark}

\begin{document}


\title{Scaling limits of nonlinear functions of random grain model, \\
with application to Burgers' equation }
\author{Donatas Surgailis}
\date{\today
\\  \small
\vskip.2cm
Vilnius University, Faculty of Mathematics and Informatics,
Naugarduko 24, 03225  Vilnius, Lithuania \\
}
\maketitle

\begin{abstract}

\medskip

We study scaling limits of nonlinear functions $G$ of random grain model $X$ on $\R^d $ with long-range dependence and marginal Poisson distribution.
Following \cite{kaj2007}  we assume that the intensity $M$ of the underlying Poisson process of grains increases together
with the scaling parameter $\lambda$ as $M = \lambda^\gamma $, for some $\gamma > 0$.
The results are applicable
to the Boolean model and exponential $G$  and rely on an expansion of $G$ in Charlier polynomials and a generalization
of Mehler's formula. Application to solution of Burgers' equation with initial aggregated random grain data is discussed.

\end{abstract}

{\small

\noi {\it Keywords:} Random grain model, nonlinear function, Boolean model, long-range dependence,  scaling limit, aggregation,
Poisson distribution, Charlier polynomials, Mehler's formula, Gaussian random field, stable random field, Burgers' equation, initial random grain data

}

\section{Introduction}

Limit theorems for random fields (RFs) with continuous or discrete $d$-dimensional argument have been extensively studied in the literature. Quite often,
such work refer to the limit distribution  of integrals
\begin{eqnarray} \label{Xlambda0}
X_{\lambda}(\phi)
&:=&\int_{\R^d}  X(\mbt) \phi(\mbt/\lambda) \d \mbt, \qquad \text{as} \  \lambda \to \infty,
\end{eqnarray}
(or respective sums in the discrete argument case), where $X = \{X(\mbt); \mbt \in \R^d \} $ is a given stationary RF,
for each $\phi$ from  a class of (test) functions $ \Phi = \{\phi: \R^d \to \R\}$. A suitably normalized
limit of \eqref{Xlambda0} is a RF  indexed by $\phi \in \Phi$, called the (isotropic) scaling limit of $X$ in this paper. The above approach
is common in the theory of generalized RFs,  and is discussed in \cite{dob1980} together with various classes of $\Phi$. 
For $\phi (\mbt) = \I (\mbt \in  A) $, where $A$ is a bounded Borel subset of $\R^d  $,
\eqref{Xlambda0} is the integral (sum) of $X(\mbt)$'s over $\mbt \in \lambda A$ whose limit distribution was studied in \cite{lah2016} for linear RF $X$ on $\Z^d $  and $A$ having irregular boundary. For $\Phi = \Phi_{{\rm rec},d} :=  \{ \I (\mbt \in ]\boldsymbol{0}, \mbx]); \mbx \in \R^d_+ \}, $
$]\boldsymbol{0}, \mbx] := \prod_{i=1}^d ]0, x_i] $,  \eqref{Xlambda0} present a $d$-dimensional analog of the partial integral
process of time series, leading to a limit RF indexed by $\mbx \in \R^d_+$.

A stationary RF $X$ on $\R^d $ with finite variance with  is said {\it long-range dependent} (LRD) if its covariance function $r_X(\mbt)
:= {\rm Cov}(X(\boldsymbol{0}), X(\mbt)) $ is not integrable: $\int_{\R^d} |r_X(\mbt)| \d \mbt = \infty $. It is well-known
that LRD RFs can  display a variety of Gaussian and non-Gaussian limit behaviors and scaling limits, see  \cite{dobmaj1979, leo1999, samo2016, pipi2017}
and the references therein. The framework in \eqref{Xlambda0} can be modified 
by replacing
the isotropic scaling $\mbt \to \mbt/\lambda $ in the test function with {\it anisotropic scaling} $\mbt \to \lambda^{-\Gamma} \mbt, $ where $\Gamma = {\rm diag}(\gamma_1, \cdots, \gamma_d), (\gamma_1, \cdots, \gamma_d) = \mbgamma \in \R^d_+ $. It was observed in \cite{ps2015, ps2016, pils2016, pils2017, sur2019, sur2020} that in dimensions $d = 2, 3$ and a rich class of LRD RFs (including Gaussian and linear RFs) the corresponding anisotropic limits
in \eqref{Xlambda0} exist for any $\mbgamma \in \R^d_+, \phi \in \Phi_{{\rm rec},d}
$ and depend on $\mbgamma $. Particularly, when $d=2 $ a {\it scaling transition} may occur meaning that there exists a critical $\gamma^0 >0$
depending on internal parameters of RF $X$
such that the limit RFs do not depend on $\mbgamma = (\gamma_1, \gamma_2 )$ for
$\gamma_2/\gamma_1 > \gamma^0$ and $\gamma_2/\gamma_1 < \gamma^0$ and exhibit a trichotomy of the (anisotropic) scaling behavior
\cite{ps2015, ps2016,  pils2017, sur2020}.

The above scaling procedures can be extended to include  aggregation, as follows. Let $X_i, i=1,2,\cdots, M $ be independent copies of
RF $X$ in \eqref{Xlambda0}. Consider the limit distribution of the `aggregate'
\begin{equation}\label{aggre1}
X_{\lambda, M}(\phi)
\ := \  \sum_{i=1}^M \int_{\R^d}  X_i(\mbt) \phi(\mbt/\lambda) \d \mbt
\end{equation}
as $\lambda \to \infty $ and $M$ increases with $\lambda $ at a certain rate.
In as follows, we refer to the limit distribution of the sums in \eqref{aggre1} 
as
the {\it scaling limit  with aggregation}. By concretising $M$
as $M = [\lambda^\gamma y]$  for some $\gamma >0, y >0 $ we see that
\eqref{aggre1} may be regarded as a discretized version of the integral in \eqref{Xlambda0}
for stationary RF $X'$ on  $\R^d \times \Z $ given by $X' (\mbt, s) :=  X_s (\mbt), (\mbt, s )  \in \R^d \times \Z $,
corresponding to (partly) anisotropic scaling $(\mbt, s) \to (\mbt/\lambda, s/\lambda^\gamma)  $ and
test function $\phi'(\mbt,s) := \phi(\mbt)  \I ( s \in ]0, y])$.
Scaling limits with aggregation in dimension $d=1$ were studied in
\cite{gaig2003, kajt2008, leip2022, miko2002, miko2007, pils2014, pilss2020, pipi2004}  and other works, in connection with applications in communications
and econometrics. As noted in \cite{leip2022}, the trichotomy of the limit distribution observed in these papers can be interpreted as a scaling transition for RFs on the plane.
Iterated limits of \eqref{aggre1}  when $M \to \infty $ (first) and $\lambda \to \infty $ (second) and/or vice versa, were discussed
in \cite{barc2017, ps2010, ps2016} and other above cited references.  A particularly simple form of aggregation occurs in models based on Poisson process
as in \eqref{RG} where summing over $M$ independent copies reduces to multiplying by $M$ the intensity of the underlying Poisson process
\cite{bier2010, bier2018, kaj2007, kajt2008, leip2022, miko2002}.


Unsurprisingly, most existing results on scaling limits of LRD RFs ($d\ge 2$) (with or without aggregation)
apply to linear models. A notable exception is Gaussian subordinated
RFs (written as a nonlinear function $G(X(\mbt))$ of a Gaussian LRD RF $X$), treated via Hermite expansion in the classical work \cite{dobmaj1979}.
For more recent developments on chaos expansions and limit theorems under Gaussian subordination, see  \cite{pec2011}.
This is in contrast to the one-dimensional
case $d = 1$, where the martingale approach developed in \cite{hoh1997} is applicable to nonlinear functions and statistics of LRD moving averages. See monographs \cite{douk2003, gir2012}. Scaling  limits of polynomial functions $G$ of LRD moving-average RFs $X$ are discussed in
\cite{sur1982, pils2017}. The case of indicator functions  $G(x) = \I(x \le y) $ for the same class of RFs was considered in \cite{douk2002}.
See \cite{kou2016} for statistical application.

A class of LRD RFs whose trajectories and distribution is very different from Gaussian and moving-average RFs are
{\it random grain} (RG) models, defined as follows.
Let $ \{ \mbu_j; j \ge 1 \} \subset \R^d$ be a Poisson point process with intensity $ \d \mbu$.
Let be given a bounded Borel set $\Xi^0 \subset \R^d$ (`generic grain')  and
and an i.i.d. sequence  $R, R_j; j \ge 1$ of r.v.s with values in $\R_+ $, distribution $F(\d r)$,
finite expectation $\E R = \int_{\R_+}  r F(\d r) < \infty$, and independent of  the Poisson process.
The RG RF  is obtained by counting at each point $\mbt$ the number of
randomly dilated grains $\mbu_j + R_j^{1/d} \Xi^0$ `centered' at $\mbu_j$ which cover it, viz.,
\begin{equation}\label{RG}
X(\mbt) := \sum_{j=1}^\infty \I (\mbt - \mbu_j \in R_j^{1/d} \Xi^0), \qquad \mbt \in  \R^d.
\end{equation}
For $d = 1$ and $\Xi^0 = ]0,1]$, \eqref{RG} is  the number of customers in a stationary M/G/$\infty$ queue with service time distribution $F (\d r)  $.
The RF in \eqref{RG}  is infinitely divisible and has a Poisson stochastic
integral representation
\begin{equation}\label{RG1}
X(\mbt) = \int_{\R^d \times \R_+} \I (\mbt -\mbu  \in r^{1/d} \Xi^0) {\cal N}(\d \mbu, \d r), \qquad \mbt \in  \R^d,
\end{equation}
where 
${\cal  N}(\d \mbu, \d r)$ is a Poisson random measure  with mean $\E {\cal N}(\d \mbu, \d r) =  \d \mbu F(\d r) $. RF in \eqref{RG1}
is  stationary and has marginal Poisson distribution with mean $ \mu := \E X(\mbt)
=\int_{\R^d}  \P (\mbt - \mbu  \in R^{1/d} \Xi^0) \d \mbu  = \Leb_d (\Xi^0) \E R $. It is well-known \cite{kaj2007} that under mild conditions on $\Xi^0$ the RG model is LRD if the distribution $F$ varies regularly at infinity with exponent $\alpha \in  (1,2)$. The conditions on $F $ and $\Xi^0$
in this paper (Assumption LRD in Sec.2)  imply that
\begin{eqnarray} \label{covRG2}
r_X (\mbt)
&\sim&\|\mbt\|^{-d(\alpha -1)} \ell(\frac{\mbt}{\|\mbt\|}),  \qquad |\mbt| \to \infty,  \quad 1 <  \alpha < 2 \nn
\end{eqnarray}
where $\ell(\mbz), \|\mbz\| =1 $ is a continuous (angular) function on the unit sphere of $\R^d $ given in \eqref{elldef} and
$\int_{\R^d} |r_X(\mbt)| \d \mbt = \infty $ holds for any $d \ge 1, \alpha \in (1,2)$.

Scaling limits with aggregation of RG model in \eqref{RG1} were discussed in
\cite{kaj2007} (see also  \cite{bier2010, bier2018}).
Let
\begin{equation} \label{Xlambda1}
X_{M}(\mbt) :=
\int_{\R^d \times \R_+}  \I(\mbt - \mbu \in r^{1/d}\,\Xi^0) {\cal N}_{M} (\d \mbu,  \d r), \qquad \mbt \in \R^d
\end{equation}
be random grain model with Poisson intensity  $
M \d \mbu F(\d r)$, a multiple of the intensity in \eqref{RG1}, where $M >0$.
\cite{kaj2007} proved that when $M$ increases with $\lambda $ at certain rate (e.g., $M = \lambda^\gamma$  for some $\gamma > 0$),
scaling limits of \eqref{Xlambda1} exhibit a trichotomy depending on $\gamma $ (described in sec. 2).

The main object of this paper are scaling limits of subordinated RFs
\begin{eqnarray}\label{Ylambda}
Y(\mbt) = G(X(\mbt)) \quad  \text{and} \quad Y_{M}(\mbt)&:=&G\big(\frac{X_{M}(\mbt) - \E X_{M}(\mbt)}{M^{1/2}} \big), \quad \mbt \in \R^d,
\end{eqnarray}
where  $X(\mbt),  X_M(\mbt)$ are as in \eqref{RG1}, \eqref{Xlambda1} and
$G(x)$ is a nonlinear function with $\E G(X(\boldsymbol{0}))^2 < \infty$ and $\E Y^2_M (\boldsymbol{0}) =
\E G^2_M(X_M(\boldsymbol{0})) < \infty,
G_M(x) := G((x- M\mu)/M^{1/2})$. 
 In other words, we discuss distributional limits of  normalized integrals
\begin{equation}\label{Ylambda1}
Y_\lambda(\phi) = \int_{\R^d} Y(\mbt) \phi(\mbt/\lambda) \d \mbt \quad \text{and} \quad
Y_{\lambda,M}(\phi) = \int_{\R^d} Y_M(\mbt) \phi(\mbt/\lambda) \d \mbt
\end{equation}
for a large class of test functions $\phi$,
as $\lambda \to \infty $ and $M = \lambda^\gamma$  for some $\gamma > 0$.  Our main result -
Theorem \ref{thmRGnon} - says that the scaling limits of \eqref{Ylambda} and \eqref{Xlambda1} are essentially the same (including the trichotomy of the
limit distribution) {\it  provided the Hermite rank of $G$ is 1}, or
\begin{equation}\label{hG}
h_{G, \mu}(1) := \E [G(Z_\mu) Z_\mu ] \ne 0, \qquad Z_\mu \sim N(0, \mu),
\end{equation}
in which case the difference between the scaling limits of \eqref{Ylambda} and \eqref{Xlambda1} reduces to the multiplicative factor in \eqref{hG}.
We also prove a similar result for scaling limits of $Y(\mbt) = G(X(\mbt))$ with \eqref{hG} replaced by the first coefficient of the expansion of $G$
in Charlier polynomials (Proposition \ref{propGlim0}).

The proofs of our results are rather simple, relying on expansion of bivariate Poisson distribution in orthogonal Charlier polynomials
and a Poissonian analog of Mehler's formula (Lemma \ref{lemmehl}). We hope that this approach can be useful to other Poisson-based RF models
and nonlinear triangular arrays. Some open problems are indicated in  Remarks \ref{remran}, \ref{remsmall}, \ref{remshot}, and \ref{remcox}.

In this paper, the limit results are applied to two subordinated RG models. The first one is $ G(x) = x \wedge 1,  \ x \in \N  $,  or
$\hat X(\mbt) := X(\mbt) \wedge 1, \mbt \in \R^d $,  referred to below as the {\it Boolean model}. By taking $\phi (\mbt) = \I(\mbt \in A) $,
where $A \subset \R^d $ is a bounded Borel set,
we see that
\begin{equation}\label{hatXA}
\int_{\R^d} \I(\mbt/\lambda \in A)    \hat X(\mbt) \d \mbt =  \Leb_d
({\cal X} \cap \lambda A) =: \hat X_\lambda (A)
\end{equation}
is the `volume' of the intersection of the {\it Boolean set}
${\cal X} := \bigcup_{j=1}^\infty (\mbu_j + R_j \Xi^0)  \subset   \R^d $
with large `inflated' set  $\lambda A$.
We remark that the Boolean model is  a basic model in stereology and stochastic geometry  \cite{sto1989}.
According to Corollary \ref{corBool},
under the above assumptions on \eqref{RG},
the limit distribution of \eqref{hatXA} is asymmetric $\alpha$-stable for a general set $A$.
The second example is the exponential function $G(x) = \e^{a x},   x \in \R $, where  $a \ne 0$ is a real parameter.
We have that $ h_{G, \mu}(1) =  a \e^{a^2 \mu /2} $ satisfies \eqref{hG} and
Theorem \ref{thmRGnon} applies to the above $G$,
see Example \ref{ex1}.  The case $G(x) = \e^{a x}$  is particular interest to the study of scaling limits
of statistical solution of Burgers' equation with random linear data as in \eqref{Xlambda1} discussed in the last
Section 4.

\noi {\it Notation.} In this paper, $\limd $ (respectively, $\limfdd$)
denote respectively the weak convergence  of distributions (respectively,
finite dimensional distributions). $C$ stands for a generic positive constant which may assume different values at various locations and whose precise value has no importance.
$\R^d := \{ \mbt = (t_1, \cdots, t_d): t_i \in \R, i=1, \cdots, d \}, \boldsymbol{0} = (0, \cdots, 0) \in \R^d,
\R^d_+ := (0,\infty)^d$.
$\I(A)$ denotes the indicator  function of a set $A$.

\section{Scaling limits of RG model}

Most results in this section either belong to \cite{kaj2007}, or are variations of the latter work. For reader's convenience, short complete proofs
are included. The following assumption is used throughout this paper without further notice.

\noi {\bf Assumption LRD} $\Xi^0 \subset \R^d$ is a bounded Borel set whereas  $F(\d r) = f(r) \d r $ has a density function such that
\begin{eqnarray}\label{falpha}
f(r)&\sim&c_f r^{-1-\alpha}, \quad r \to \infty  \qquad (\exists \, c_f >0, \ \alpha \in  (1,2)).
\end{eqnarray}
Moreover, the function  $ (r, \mbz) \mapsto \Leb_d \big(\Xi^0 \cap (\Xi^0 - r^{-1/d} \, \mbz)\big) $ is continuous in  $ (r, \mbz) \in \R_+  \times \{ \|\mbz \| = 1 \} $ and nontrivial.

\smallskip

\noi
The class of test functions in \eqref{Xlambda0} and elsewhere  in this paper is
\begin{equation}\label{Phi}
\Phi   =   L^1 (\R^d) \cap L^\infty (\R^d),
\end{equation}
where $L^1(\R^d) $ ($L^\infty (\R^d)$)
stand for the linear space of all Borel functions $\phi: \R^d \to \R$ such that
$\int_{\R^d} |\phi(\mbt)| \d \mbt < \infty $ (respectively, such that $\sup_{\mbt \in \R^d} |\phi(\mbt)| < \infty$). The proofs of
Theorems \ref{thmRG} and \ref{thmRGnon} use the fact 
that any Borel function is a.e. continuous on $\R^d$ (Lusin's theorem).

\begin{proposition} \label{propcov}

\noi (i) Relation  \eqref{covRG2} holds with bounded, continuous and nonnegative angular function
\begin{eqnarray}  \label{elldef}
\ell(\mbz)&:=&c_f \int_{\R_+} \Leb_d \big(\Xi^0 \cap (\Xi^0 - r^{-1/d} \, \mbz)\big) r^{-\alpha} \d r.
\end{eqnarray}
\noi (ii) For any $\phi \in \Phi $  as $\lambda \to \infty $
\begin{eqnarray}\label{c1}
c_\lambda(\phi) :=   \int_{\R^{2d}} \phi(\mbt_1/\lambda) \phi(\mbt_2/\lambda) r_X(\mbt_1 - \mbt_2) \d \mbt_1 \d \mbt_2
&\sim&
\lambda^{d(3 - \alpha)} c(\phi),
\end{eqnarray}
where
\begin{equation}  \label{c2}
c(\phi)  := \int_{\R^{2d}} \phi(\mbt_1) \phi(\mbt_2) \ell \big(\frac{\mbt_1 - \mbt_2}{\|\mbt_1 -  \mbt_2\|}\big) \frac{\d \mbt_1 \d  \mbt_2}
{\|\mbt_1 -  \mbt_2\|^{d(\alpha -1)}}
\end{equation}
and the integral on the r.h.s. of \eqref{c2} converges.
Moreover,
\begin{eqnarray}\label{c3}
\int_{\R^{2d}} \phi(\mbt_1/\lambda) \phi(\mbt_2/\lambda) (1 \wedge  r^2_X(\mbt_1 - \mbt_2)) \d \mbt_1 \d \mbt_2
&=&\begin{cases}
O(\lambda^d), &\alpha > 3/2, \\
O(\lambda^{2d(2-\alpha)}), &\alpha < 3/2, \\
O(\lambda^d \log \lambda), &\alpha = 3/2.
\end{cases}
\end{eqnarray}

\end{proposition}

\noi {\bf Proof.} (i) From \eqref{Xlambda1} we have that
\begin{eqnarray} \label{covRG}
\|\mbt \|^{d(\alpha -1)} r_X(\mbt)
&=&\|\mbt \|^{d(\alpha -1)}
\int_0^\infty \Leb_d (r^{1/d} \Xi^0 \cap (r^{1/d}  \Xi^0 - \mbt)) f(r) \d r    \\
&=&\int_0^\infty  \Leb_d \big( \Xi^0 \cap  (\Xi^0 - r^{-1/d} \frac{\mbt}{\|\mbt\|})\big)  \tilde f(r; \|\mbt\|^d)   r^{-\alpha} \d r \nn
\end{eqnarray}
where, for any $ r>0$,  according to \eqref{falpha}
\begin{equation}\label{tildef}
\tilde f(r; \lambda ) := (\lambda r)^{1+ \alpha} f(\lambda r) \ \to \  c_f  \qquad (\lambda \to  \infty).
\end{equation}
By boundedness of $\Xi^0$, $\Leb_d \big(\Xi^0 \cap (\Xi^0 - r^{-1/d} \, \mbz)\big) \le \Leb_d (\Xi^0) $ is
bounded and vanishes for $r >0$ small enough
uniformly in $\|\mbz \| = 1 $. Thus, \eqref{covRG2}  follows from \eqref{covRG} and \eqref{tildef} by the dominated convergence theorem. Properties
of $\ell (\mbz) $ follow easily  from its definition and Assumption LRD.

\smallskip

\noi (ii)  Following \eqref{covRG}, write
\begin{eqnarray}\label{c4}
\frac{c_\lambda(\phi)}{\lambda^{d(3-\alpha)}}
&=&\int_{\R^{2d}} \phi(\mbt_1) \phi(\mbt_2) \|\mbt_1 - \mbt_2\|^{d(1-\alpha)} \d \mbt_1 \d \mbt_2 \nn \\
&& \times \int_0^\infty
 \Leb_d \Big(\Xi^0 \cap \big(\Xi^0 - r^{-1/d} \, \frac{\mbt_1 - \mbt_2}{\|\mbt_1 - \mbt_2\|}\big)\Big) \tilde f(r; \lambda \|\mbt_1 - \mbt_2\|)   r^{-\alpha} \d r \nn
\end{eqnarray}
and use \eqref{tildef} and the argument as in the proof of (i) to conclude \eqref{c1}.

The proof of \eqref{c3} is  similar but simpler. Indeed, by boundedness of $\Xi^0 $ and \eqref{falpha} we have that
$|r_X(\mbt)| \le  C(1 \wedge \|\mbt\|^{d(1-\alpha)}), $ c.f.  \eqref{covRG}, so  that
the l.h.s.  of \eqref{c3} does not  exceed $C \int_{\R^{2d}} |\phi(\mbt_1/\lambda) \phi(\mbt_2/\lambda)| (1 \wedge \|\mbt_1 - \mbt_2 \|^{2d(1  - \alpha)} )
\d \mbt_1 \d \mbt_2 $ whose evaluation by the r.h.s. in \eqref{c3} for $\phi \in \Phi $ in \eqref{Phi}
is elementary.  \hfill $\Box$

\smallskip

Introduce a Gaussian RF $B_\alpha(\phi)$ indexed by functions $\phi  \in  \Phi $ as stochastic integral
\begin{eqnarray} \label{BB}
B_\alpha(\phi)
&:=&\int_{\R_+ \times \R^d}  W_\alpha(\d r, \d \mbu) \int_{\R^d} \phi (\mbt) \I ( \mbt - \mbu \in r^{1/d} \Xi^0) \d \mbt, \qquad
\phi \in \Phi
\end{eqnarray}
w.r.t. Gaussian white noise $W_\alpha$  with zero mean and variance
$\E W_\alpha(\d r, \d \mbu)^2 = c_{f} r^{-1-\alpha} \d r \d \mbu  $.
Observe that the variance
\begin{eqnarray}
\E B_\alpha(\phi)^2
&=&c_{f} \int_{\R_+ \times \R^d} r^{-1-\alpha} \d r \d \mbu \Big(\int_{\R^d} \phi (\mbt) \I ( \mbt - \mbu \in r^{1/d} \Xi^0)
\d \mbt \Big)^2 \nn \\
&=&c_{f}\int_{\R^{2d}} \phi (\mbt) \phi(\mbt') \d \mbt \d \mbt'
\int_{\R_+ \times \R^d} \I ( \mbt - \mbu \in r^{1/d} \Xi^0,  \mbt' - \mbu \in r^{1/d} \Xi^0 )    r^{-1-\alpha} \d r \d \mbu \nn \\
&=&c(\phi) \label{BBvar}
\end{eqnarray}
coincides with \eqref{c2}.  Let be given a centered
Poisson  random  measure
$\widetilde N_\alpha(\d r, \d \mbu) = N_\alpha(\d r, \d \mbu) - \E N_\alpha(\d r, \d \mbu)$ on $\R^d \times \R_+$
with variance $\E \widetilde N_\alpha(\d r, \d \mbu)^2 = \E N_\alpha(\d r, \d \mbu) =
c_{f} r^{-1-\alpha} \d r \d \mbu  $  the same as the variance of $W_\alpha$, and an $\alpha$-stable random measure  $L_\alpha (\d \mbu) $
on $\R^d $ with the characteristic function
\begin{eqnarray} \label{Lalpha}
\E \e^{\i \theta L_\alpha (A)}
&:=&\exp \big\{ \Leb_d (A) \, c_f \int_0^\infty (\e^{\i \theta r}  - 1 - \i \theta r) r^{-1-\alpha} \d r \big\} \\
&=&\exp \big\{- \Leb_d (A) \, \sigma_\alpha |\theta|^\alpha \big(1 - \i \, {\rm sgn}(\theta)\, \tan (\pi \alpha/2)\big)\big\}, \quad \theta \in \R, \nn
\end{eqnarray}
where $A \subset \R^d $ is any Borel set with $\Leb_d (A)< \infty $ and $\sigma_\alpha :=  ...  $.
Introduce RFs $L_\alpha  $ and $J_\alpha $ as stochastic integrals
\begin{eqnarray}\label{JJ}
L_{\alpha}(\phi)
&:=&\int_{\R^d}  \phi (\mbu) L_{\alpha}(\d \mbu), \\
J_\alpha(\phi)
&:=&\int_{\R_+ \times \R^d} \widetilde N_\alpha (\d r, \d \mbu) \int_{\R^d} \phi (\mbt) \I ( \mbt - \mbu \in r^{1/d} \Xi^0) \d \mbt \nn
\end{eqnarray}
w.r.t. to the above random measures which are well-defined for any $\phi \in \Phi$.
Clearly, $\E J_\alpha(\phi)^2 = \E B_\alpha(\phi)^2 $. Denote $\Psi (z) := \e^{\i z} - 1 - \i z, z \in  \R$.

\begin{theorem} \label{thmRG} Let $ X_{\lambda, M}(\phi)= \int_{\R^d}  X_M(\mbt) \phi(\mbt/\lambda) \d \mbt$ with $X_M $   as in \eqref{Xlambda1}, $M  = \lambda^{\gamma} \ (\gamma > 0)$. Then
for any $\phi \in  \Phi $  as $\lambda \to \infty $
\begin{eqnarray}\label{limXinfty}
\lambda^{-H(\gamma)}(X_{\lambda, M}(\phi) -
\E X_{\lambda, M}(\phi)) &\limd&
\begin{cases}
B_\alpha(\phi),  &  \gamma > d(\alpha-1), \ H(\gamma) = \frac{\gamma +  (3-\alpha)d}{2},  \\
L_{\alpha}(\phi), &   \gamma < d(\alpha-1),  \   H(\gamma) = \frac{\gamma + d}{\alpha},  \\
J_\alpha(\phi),  & \gamma = d(\alpha-1), \ H(\gamma) = d.
\end{cases}
\end{eqnarray}

\end{theorem}

\noi  {\bf Proof.} Let $j_\lambda (\theta) :=   \log \E \exp \{ \i  \theta \lambda^{-H(\gamma)}
 (X_{\lambda, M}(\phi) - \E X_{\lambda, M}(\phi))  \}, \theta \in \R$  denote the log-characteristic functional
of the l.h.s. in  \eqref{limXinfty}.  We need to show that it converges to the corresponding functional $j(\theta)$   of  the r.h.s. as
$\lambda \to \infty $.

\noi \underline{Case $\gamma > d (\alpha-1)$.} Set
$\tilde j_\lambda (\theta) :=    - (\theta^2/2) {\rm Var} (\lambda^{-H(\gamma)} X_{\lambda, M }(\phi)), $  $
j (\theta)
:=  \log \E \exp \{ \i  \theta
B_\alpha(\phi)\} = -
(\theta^2/2) \E B_\alpha(\phi)^2 $. Then $\tilde j_\lambda (\theta) \to  j(\theta) \ (\lambda \to \infty)$, see \eqref{c1}, \eqref{c2}, \eqref{BBvar}.
Next,
\begin{eqnarray*}
j_\lambda (\theta) - \tilde j_\lambda(\theta)
&=&\lambda^{\gamma}
\int_{\R^d \times \R_+} \d \mbu f(r) \d r
\Big\{\Psi \Big(\frac{\theta}{\lambda^{H(\gamma)}} \int_{\R^d} \phi(\mbt/\lambda) \I (\mbt - \mbu \in r^{1/d} \Xi^0) \d \mbt \Big) \\
&&\hskip4cm + (1/2) \Big(\frac{\theta} {\lambda^{H(\gamma)}}  \int_{\R^d} \phi(\mbt/\lambda)
\I (\mbt - \mbu \in r^{1/d} \Xi^0) \d \mbt \Big)^2 \Big\}
\end{eqnarray*}
Using $|\Psi (z) + (1/2) z^2 | \le |z|^3 $ and \eqref{c1} we get
\begin{eqnarray*}
|j_\lambda (\theta) - \tilde j_\lambda(\theta)|
&\le&C\lambda^{\gamma} \int_{\R^d \times \R_+}
\Big|\lambda^{-H(\gamma)} \int_{\R^d} \phi(\mbt/\lambda) \I (\mbt - \mbu \in r^{1/d} \Xi^0) \d \mbt \Big|^3  \d \mbu f(r) \d r  \\
&\le&C\lambda^{-(\gamma/2) - (3/2)(3-\alpha)d}  \int_{\R^d} |\phi(\mbt/\lambda)| \d \mbt \, {\rm Var} (X_{M}(\phi)) \\
&\le&C \lambda^{-(\gamma - d(\alpha-1))/2} \ = \ o(1)
\end{eqnarray*}
since  $|\int_{\R^d} \phi(\mbt/\lambda) \I (\mbt - \mbu \in r^{1/d} \Xi^0) \d \mbt|
\le \int_{\R^d} |\phi(\mbt/\lambda)| \d \mbt < C \lambda^d $ uniformly in $\mbu \in \R^d$.  Thus,
$j_\lambda (\theta) \to  j(\theta) \, \forall \theta \in \R$.

\smallskip

\noi \underline{Case $\gamma < d (\alpha-1)$.} We have
\begin{eqnarray} \label{jrhoeq}
j_\lambda (\theta)&=&\int_{\R^d \times \R_+} \Psi \Big(
\theta \lambda' \int_{\R^d} \phi(\mbt + \mbu)
\I (\mbt \in (r/\lambda')^{\frac{1}{d}}\, \Xi^0) \d \mbt \Big) \tilde f(r; \lambda^{\frac{\gamma + d}{\alpha}})
\frac{ \d \mbu \d r}{r^{\alpha +1}},
\end{eqnarray}
where
$\lambda' := \lambda^{d - \frac{\gamma +d}{\alpha}} \to  \infty  $ and $\tilde f $ is as in \eqref{tildef}.
Since $\phi$ is a.e. continuous and 
$h_\lambda (\mbu, r) := \lambda' \int_{\R^d} \phi(\mbt + \mbu) \I (\mbt  \in (\frac{r}{\lambda'})^{\frac{1}{d}}\, \Xi^0)  \d \mbt \to r \phi(\mbu) \Leb_d (\Xi^0) =: h(\mbu, r) $ at each continuity point $\mbu $ of $\phi$, we get that 
$h_\lambda (\mbu, r)\to h(\mbu, r) $ 
for a.e. $\mbu \in \R^d $ and
each $r >0$. This fact together with  \eqref{jrhoeq} and \eqref{tildef} suggest that
\begin{eqnarray}\label{jrholim}
j_\lambda (\theta)
\ \to\ j(\theta)
&=&c_f \int_{\R^d \times \R_+} \big(\e^{\i \theta r \phi(\mbu) \Leb_d (\Xi^0)} - 1 - \i \theta r \phi(\mbu) \Leb_d (\Xi^0)\big)
\frac{\d \mbu \d r}{r^{\alpha+1}} \\
&=&\log \E \exp\{\i \theta L_{\alpha} (\phi)\}, \nn
\end{eqnarray}
see \eqref{Lalpha}, \eqref{JJ}. The convergence in \eqref{jrholim} can be justified using Pratt's lemma \cite{pra1960}, as follows.
It suffices to show that
\begin{eqnarray}
\int_{\R^d} \Psi ( h_\lambda (\mbu, r)) \d \mbu
&\to & \int_{\R^d} \Psi ( h(\mbu, r))  \d \mbu,  \qquad \forall  \,  r>0,  \label{pr1} \\
\int_{\R^d} |\Psi ( h_\lambda (\mbu, r))| \d \mbu &\le &C (r \wedge r^2), \label{pr2} \\
\int_0^\infty  (r \wedge r^2) \tilde f(r; \lambda^{\frac{\gamma + d}{\alpha}}) \frac{ \d r}{r^{\alpha +1}}
&\to& c_f \int_0^\infty  (r \wedge r^2)\frac{\d r}{r^{\alpha +1}}. \label{pr3}
\end{eqnarray}
Relations \eqref{pr1} and \eqref{pr3} are rather easy.  To show \eqref{pr2}, use
$|\Psi(z)| \le 2(|z|  \wedge |z|^2)$ and boundedness of $\Xi^0$. Thus, the l.h.s. of \eqref{pr2} does not exceed 
$  C \int_{\R^d} \min
\Big(\lambda' \int_{\R^d} |\phi(\mbt + \mbu)| \I (\|\mbt \| \le  C (\frac{r}{\lambda'})^{\frac{1}{d}})  \d \mbt,
\big(\lambda' \int_{\R^d} |\phi(\mbt + \mbu)| \I (\|\mbt \| \le  C (\frac{r}{\lambda'})^{\frac{1}{d}})  \d \mbt\big)^2 \Big) \d \mbu $
$\le C \min \Big( \lambda' \int_{\R^d} (\int_{\R^d} |\phi(\mbt + \mbu)| \d \mbu ) \I (\|\mbt \| \le  C (\frac{r}{\lambda'})^{\frac{1}{d}})  \d \mbt,
(\lambda')^2 \int_{\R^{2d}} (\int_{\R^d} |\phi(\mbt_1 + \mbu)  \phi (\mbt_2 + \mbu)|  \d \mbu )
\I (\|\mbt_1 \| \le  C (\frac{r}{\lambda'})^{\frac{1}{d}}, \|\mbt_2 \| \le  C (\frac{r}{\lambda'})^{\frac{1}{d}}  )  \d \mbt_1 \d \mbt_2 \Big)
$
$\le C \|\phi\|_{L^1} \big(\lambda' \int_{\R^d}  \I (\|\mbt \| \le  C (\frac{r}{\lambda'})^{\frac{1}{d}})  \d \mbt\big) \wedge
\big(\lambda' \int_{\R^d}  \I (\|\mbt \| \le  C (\frac{r}{\lambda'})^{\frac{1}{d}})  \d \mbt\big)^2 \le
C (r \wedge r^2) $. This proves \eqref{pr2} and \eqref{jrholim}.

\smallskip

\noi \underline{Case $\gamma = d (\alpha-1)$.} The expression in  \eqref{jrhoeq} with $\lambda' =1$ for $j_\lambda (\theta)$
is valid. The proof of \eqref{limXinfty} in this case is similar to that when  $\gamma < d (\alpha-1)$
and we omit the details. \hfill $\Box$

\begin{rem}
{\rm Theorem \ref{thmRG} applies also to $\gamma =0 $ or integrals in \eqref{Xlambda0}
of the RG model in \eqref{RG1} with Poisson intensity  $\d \mbu F(\d r)$. Indeed, the argument in the case $\gamma < d (\alpha-1)$ applies without change when  $\gamma = 0$, yielding the stable limit
\begin{eqnarray}\label{limX0}
\lambda^{-d/\alpha}(X_{\lambda}(\phi) -
\E X_{\lambda}(\phi)) &\limd&
L_{\alpha}(\phi), \qquad \forall \, \phi \in \Phi.
\end{eqnarray}
}
\end{rem}

\begin{rem}
{\rm Condition $M = \lambda^\gamma $ in Theorem \ref{thmRG} can be weakened \cite{kaj2007}. Particularly, it can be replaced by $M = y \lambda^\gamma $ for some $y>0$, in which case the  control measure of the  limit RFs in  \eqref{limXinfty} contain the extra multiplicative factor $y$.
For $d=1, \phi_s(t)  := \I(t \in ]0,s]), t\in \R, (s,y) \in \R^2_+$ we have that $X_{\lambda, M} (\phi_s)
= \int_0^{\lambda s} X_M (t) \d t $, where $X_M (t) = \sum_{j \in \Z} \I( u_j (M) <  t \le  u_j(M) + R_j) $  is the number of customers at time $t$
in the M/G/$\infty$   queue with service time  distribution $R $  and a Poisson arrival process $\{u_j (M); j \in \Z \}$ with  intensity $M \d u
= y \lambda^\gamma \d u $. Then
 \begin{eqnarray}\label{limXy}
\lambda^{-H(\gamma)}\int_0^{\lambda s} (X_{M} (t) - \E X_{M} (t))  \d t
&\limfdd&
\begin{cases}
B_\alpha(s,y),  &  \gamma > \alpha-1, \ H(\gamma) = \frac{\gamma +  3-\alpha}{2},  \\
L_{\alpha}(s,y), &   \gamma < \alpha-1,  \   H(\gamma) = \frac{\gamma + 1}{\alpha},  \\
J_\alpha(s,y),  & \gamma = \alpha-1, \ H(\gamma) = 1,
\end{cases}
\end{eqnarray}
where the limits are RFs indexed by $(s,y)\in \R^2_+$ written as stochastic integrals
w.r.t. Gaussian, $\alpha$-stable and Poisson random measures on $\R \times \R^2_+ $ analogously to \eqref{BB}, \eqref{JJ}. See  \cite{leip2022} for
details. Particularly,
  $\{B_\alpha(s,v); (s,v) \in \R^2_+\} $ is a fractional Brownian sheet with Hurst parameters $(H_1,H_2) = (\frac{3-\alpha}{2}, \frac{1}{2})$,
  $\{L_\alpha(s,v); (s,v) \in \R^2_+\} $ is an $\alpha$-stable L\'evy sheet, and     $\{J_\alpha(s,v); (s,v) \in \R^2_+\} $ is the Telecom RF defined
  in \cite{leip2022} as a RF extension of the corresponding Telecom process in \cite{kajt2008}.  As noted above, Theorem \ref{thmRG} is
essentially due to \cite{kaj2007} while \eqref{limXy} is a version of the results in
\cite{miko2002, kajt2008} and other previous work.

   }
\end{rem}

\section{Charlier polynomials and Mehler's formula}

The derivations in this sec. can be compared to the discussion of Hermite polynomials and expansions in the case of Gaussian distribution
\cite[pp.22-26]{gir2012}. See \cite{erd1953} for classical Mehler's  formula for Hermite polynomials.

Recall from \cite{sur1984, pec2011}
the definition of Charlier polynomials $P_k (x; \mu) $ of discrete variable $x \in \N$ through
the generating function:
\begin{equation}\label{jgen}
{\cal P}(u; x, \mu) := \sum_{k=0}^\infty \frac{u^k}{k!} P_k(x; \mu) = (1+ u)^x \e^{-u \mu},
\end{equation}
where the series is convergent for any $x \in \N, \mu  >0 $ and any (complex)  $u \in \C$. We have $P_0(x; \mu) = 1, P_1(x; \mu) = x - \mu, P_2 (x; \mu)
= x^2 -  (2\mu  + 1) x + \mu^2 $ and
\begin{equation}\label{jdef}
P_k (x; \mu) = (-1)^k \mu^k p(x; \mu)^{-1} D_-^k p(x; \mu), \qquad k \in \N
\end{equation}
where $D_-^k := D_- D_-^{k-1} $ is the backward difference operator,
$D_- G(x) := G(x) - G(x-1)\I(x \ge 1), D^0_- G(x) = G(x)$ and
\begin{equation}
p(x; \mu) = \e^{-\mu} \frac{\mu^x}{x!}, \quad x \in \N
\end{equation}
is the distribution  of Poisson r.v.  $N$ with mean $\mu$. Relation \eqref{jdef} follows from
\eqref{jgen} using the identity $(1+ u) \partial {\cal P}(u; x,\mu)/\partial u = (x- (1+u)\mu) {\cal P}(u;  x,\mu) $ \cite{sur1984}.
We have
\begin{eqnarray}\label{jorth}
&&\E P_k(N; \mu) = 0, \qquad \E P_k(N;\mu)^2 = k! \mu^k,   \quad k = 1,2, \cdots, \\
&& \E P_k (N; \mu) P_\ell (N; \mu) = 0, \quad k \neq \ell = 0,1, \cdots. \nn
\end{eqnarray}
Facts \eqref{jorth} follow from multiplying the series in \eqref{jgen} at the points $u $ and $v$ and taking the expectation of the product:
\begin{eqnarray*}
\sum_{k, \ell =0}^\infty \frac{u^k v^\ell}{k! \ell!} \E P_k (N; \mu) P_\ell (N;\mu)
&=&\e^{-(u+v)\mu} \E [((1 + u)(1+v))^N ] \\
&=&\e^{ \mu u v}\ = \ \sum_{k=0}^\infty \frac{(\mu u  v)^k}{k!} \nn
\end{eqnarray*}
and equating the coefficients of $u^k v^\ell, k, \ell  \in \N $  of the power series on both sides.

Any $G = G(x), x \in \N $ with $\E G^2 (N) < \infty $
can be uniquely  expanded in Charlier polynomials
\begin{equation}\label{Gexp}
G(x) = \sum_{k=0}^\infty \frac{c_G(k;\mu)}{k!}  P_k(x; \mu), \qquad x \in \N
\end{equation}
where
\begin{equation}\label{CG1}
c_G(k;\mu) := \mu^{-k} \E G(N) P_k(N; \mu),  \quad k \in \N
\end{equation}
are {\it  Charlier coefficients}   of $G$ in \eqref{Gexp}.
\eqref{jdef} and summation by  parts yields another expression for these coefficients
\begin{equation}\label{CG2}
c_G(k;\mu) = \E  D_+^k G(N),  \quad k \in \N,
\end{equation}
where $D^k_+ := D_+ D_+^{k-1} $ is the forward difference operator,
$D_+ G(x) := G(x+1) - G(x), D^2_+ G(x) = D_+ G(x+1) -  D_+ G(x) = G(x+2) - 2 G(x+1) + G(x)$ etc. \eqref{CG1} and \eqref{jorth} yield the bound
\begin{equation}\label{CGbdd}
|c_G(k;\mu)| \le \mu^{-k}\sqrt{ \E [G^2(N)] \E [P^2_k(N; \mu)]} = C (k!/\mu^k)^{1/2},  \quad C = \sqrt{ \E [G(N)^2]}.
\end{equation}

\begin{lemma} \label{lemmehl} Let $M_i, i=1,2,3, $ be independent Poisson distributed r.v.s
with respective means $\E M_1 = \mu_1 - \mu_3, \, \E M_2 = \mu_2 - \mu_3, \, \E M_3 = \mu_3, $  $0 \le \mu_3 < \mu_1 \wedge \mu_2 $,
 and
\begin{equation*} \label{N12}
N_{i} := M_i + M_3, \quad  i=1,2.
\end{equation*}
Let
\begin{equation}\label{pxy}
p(x, y; \mu_1,\mu_2,\mu_3) := \P(N_{1} = x, N_{2} = y), \qquad (x,y) \in \N^2
\end{equation}
denote the joint distribution of $(N_{1},  N_{2})$. 

\noi (i) {\rm (Orthogonality property):} For any $k, \ell \in \N$
\begin{eqnarray} \label{Porth}
\E P_k(N_{1}; \mu_1) P_\ell (N_{2}; \mu_2) = \begin{cases}0, &k\neq \ell, \\
\mu^k_3 k!, &k = \ell,
\end{cases}
\end{eqnarray}
with the convention $0^0 := 1 $.

\noi (ii) Let $G_i = G_i(x), x \in \N, i=1,2 $ be given functions such that
\begin{equation}\label{Gvar}
\E G^2_i (N_{i}) < \infty,  \qquad i=1,2.
\end{equation}
Then
\begin{eqnarray}
\label{Mcov2}
\E G_1(N_{1}) G_2(N_{2})
&=&\sum_{k=0}^\infty \frac{c_{G_1}(k;\mu_1) c_{G_2}(k;\mu_2) }{k!} \mu_3^k.
\end{eqnarray}

\noi (iii)  {\rm (Mehler's formula):}
\begin{eqnarray}\label{Pmehl}
p(x,y; \mu_1,\mu_2,\mu_3)
&=&\sum_{k=0}^\infty \frac{\mu^k_3}{k!} D^k_- p(x; \mu_1)  D^k_- p(y; \mu_2) \\
&=&p(x; \mu_1) p(y; \mu_2)  \sum_{k=0}^\infty \frac{\rho_{12}^k }{k!} 
P_k(x; \mu_1) P_k (y; \mu_2), \nn
\end{eqnarray}
where $\rho_{12} := \mu_3/\sqrt{\mu_1 \mu_2}  = {\rm Corr}(N_1,N_2)$ is the correlation coefficient.

\end{lemma}

\noi {\bf Proof.}  (i) The proof of \eqref{Porth} using  the generating function in \eqref{jdef}  is similar as in
the univariate case of \eqref{jorth}.
Consider the expectation
\begin{eqnarray} \label{bigen}
\E {\cal  P}(u; N_{1}, \mu_1)  {\cal P}(v; N_{2}, \mu_2)
&=&\e^{-u \mu_1- v \mu_2}  \E [(1+ u)^{N_{1}} (1+  v)^{N_{2}}]  \\
&=&\e^{-u\mu_1 - v \mu_2}  \E [(1+ u)^{M_1}] \E[ (1+  v)^{M_2}] \E [((1+u)(1+v))^{M_3}] \nn \\
&=&\e^{-u\mu_1 - v \mu_2} \e^{(\mu_1  - \mu_3)u} \e^{(\mu_2  - \mu_3) v} \e^{((1+u)(1+v) -1)\mu_3} \nn \\
&=&\e^{uv \mu_3} = \sum_{k=0}^\infty \frac{(u v \mu_3)^k}{k!}. \nn
\end{eqnarray}
On the other hand,
\begin{eqnarray*}
\E {\cal P}(u; N_{1}, \mu_1)  {\cal P}(v; N_{2}, \mu_2)
&=&\sum_{k, \ell =0}^\infty \frac{u^k v^\ell}{k!\ell! } \E [P_k(N_{1}; \mu_1) P_\ell (N_{2}; \mu_2)]
\end{eqnarray*}
and \eqref{Porth} follows  by  equating the coefficients of $u^k v^\ell, k, \ell  \in \N $  of the power series on both sides.

\smallskip

\noi (ii) Immediate from \eqref{Gexp} and \eqref{Porth}.

\smallskip

\noi (iii) Apply  \eqref{Mcov2} to $G_1 (x) :=  \I(x = n),   G_2 (x) :=  \I(x = m), $ for given $n, m \in \N $. By
\eqref{CG2}, \eqref{jdef},
$c_{G_1}(k ;\mu_1) =  \E [D^k_+ \I(N_{1} = n)] = D^k_- p(n; \mu_1) = (-1)^k \mu_1^{-k} P_k (n; \mu_1) p(n; \mu_1), $
$c_{G_2}(k;\mu_2) =  \E [D^k_+ \I(N_{2} = m)] = D^k_- p(m; \mu_2) = (-1)^k \mu_2^{-k} P_k (m; \mu_2) p(m; \mu_2). $
This proves \eqref{Pmehl} and the lemma, too. \hfill $\Box$

\begin{rem} {\rm  Let $N = \{N_t; t =0,1,\cdots \} $ be a stationary Markov process on $\N$ with marginal Poisson distribution $\P(N_t = x) = p(x; \mu)$ and
transition probabilities
\begin{equation}\label{INAR}
p(y|x; \mu) := \frac{p(x,y; \mu)}{p(x; \mu)}, \qquad x,y  \in \N,
\end{equation}
where $p(x,y;\mu) := p(x,y; \mu, \mu, \mu_3) $ is the joint distribution in \eqref{pxy} with $\mu_1 = \mu_2 =: \mu > \mu_3$.
This process is well-known in the literature as the {\it Poisson AR(1)} or {\it INAR(1)} and is related
to M/M/$\infty$ queueing system, see e.g. \cite{mcKe2003}. Substitution of  \eqref{Pmehl} into \eqref{INAR} yields an
expansion of \eqref{INAR} into a series of Charlier polynomials.
Since the transition probability of the Poisson INAR(1)  process is usually written
via a different expansion,  the coincidence  with \eqref{INAR} most easily can be verified through the bivariate generating
function of $(N_1,N_2)$ in \eqref{bigen}, see \cite[(9)]{mcKe2003}.  We remark that the Poisson INAR(1) process is a  particular case of stationary
Markov evolutions of non-interacting particle systems with Poisson marginal distribution discussed in \cite{sur1984}, closely related
to chaos expansions in multiple Poisson stochastic integrals.

}
\end{rem}

Given a $G(x), x \in \N, \E G(N)^2 < \infty $ with Charlier expansion
in \eqref{Gexp} we define the {\it Charlier rank $k^*(G;\mu)$ of $G$ } as the minimal $k\ge 1 $ such that $c_G(k;\mu) \ne 0$, viz.,
\begin{equation*}
k^*(G;\mu) := \min \{k  \ge 1:  c_G(k;\mu) \ne 0 \}.
\end{equation*}
Lemma  \ref{lemmehl} \eqref{Mcov2} and the bound in \eqref{CGbdd} imply the following

\begin{corollary} \label{corM} Let $G_i, N_i, i=1,2, $ be as in Lemma \ref{lemmehl},
$\rho_{12} = {\rm Corr}(N_1,N_2)$,
 $k^* :=  k_C(G_1;\mu_1) \vee k_C(G_2;\mu_2) $.
Then
\begin{eqnarray*}
{\rm Cov}(G_1(N_1), G_2(N_2))
&=&\sum_{k=k^*
}^\infty \frac{c_{G_1}(k;\mu_1) c_{G_2}(k;\mu_2) }{k!} \mu_3^k \nn \\
&=&\frac{c_{G_1}(k^*;\mu_1) c_{G_2}(k^*;\mu_2) }{k^*!} \mu_3^{k^*} + R(k^*), \nn
\end{eqnarray*}
where 
\begin{equation*}
|R(k^*)|\ \le \ \frac{ (\mu_3/\sqrt{\mu_1 \mu_2})^{k^* + 1}}{1 - (\mu_3/\sqrt{\mu_1\mu_2})} \prod_{i=1}^2 \E^{1/2}[ G(N_i)^2].
\end{equation*}
Moreover,
\begin{eqnarray*}
|{\rm Cov}(G_1(N_1), G_2(N_2))|
&\le&\sum_{k=k^*}^\infty |\rho_{12}|^k  \frac{|c_{G_1}(k;\mu_1)c_{G_2}(k;\mu_2)| (\mu_1 \mu_2)^{k/2} }{k!}  \\
&\le&|\rho_{12}|^{k^*} \prod_{i=1}^2  \big(\sum_{k=k^*}^\infty \frac{c^2_{G_i}(k;\mu_i) \mu_i^k }{k!} \big)^{1/2} \nn \\
&\le&|\rho_{12}|^{k^*}\sqrt{{\rm Var}(G_1(N_1)) {\rm Var}(G_2(N_2))}. \nn
\end{eqnarray*}
Particularly, $\sup \{ |{\rm Cov}(G_1(N_1), G_2(N_2))|:  {\rm Var}(G_i(N_i)) = 1, i=1, 2\} = |\rho_{12}| $
and the last supremum is attained
by linear functions $G_i(x) = x/\sqrt{\mu_i}, i=1,2 $.

\end{corollary}

\section{Scaling limits of nonlinear functions of RG model}

We study the limit distribution of
$Y_{M}(\phi), Y_{\lambda,M}(\phi)  $ defined in \eqref{Ylambda}-\eqref{Ylambda1}  as $\lambda \to \infty $ and $M = \lambda^\gamma \to \infty $, for some
$\gamma >0$, 
where $G = G(x), x \in \R $ is a general function satisfying some conditions.
We will show that the limit of  $ Y_{\lambda,M}(\phi) $ is
the same as that of the linear integral $X_{\lambda,M}(\phi)$ in Theorem  \ref{thmRG}  up to the multiplicative constant
equal to the first coefficient in the Hermite expansion of $G$.
Since for each $\mbt \in \R^d $, the quantity inside $G$ in $Y_M(\mbt)$ of \eqref{Ylambda}, viz.,
$(X_M(\mbt) - \E X_M(\mbt))/M^{1/2} \limd Z_\mu \sim N(0,  \mu) $
when $M \to  \infty $, we consider the corresponding Hermite  expansion
\begin{equation}\label{Ghermite}
G(x) =  \sum_{k=0}^\infty \frac{h_{G,\mu}(k)}{k!}  H_k(x; \mu)
\end{equation}
in Hermite polynomials $H_k (x; \mu), k\in \N, x \in \R $ with generating
function $\sum_{k=0}^\infty (u^k/k!) H_k (x; \mu) = \e^{ u x - \mu u^2/2} $ and Hermite coefficients
\begin{eqnarray}\label{Gh}
h_{G,\mu} (k) := \mu^{-k} \E [ G(Z_\mu) H_k(Z_\mu; \mu)], \qquad k \in \N.
\end{eqnarray}

\begin{theorem} \label{thmRGnon} Let $ Y_{\lambda, M}(\phi) $ be as in \eqref{Ylambda}, where $X_{M}(\mbt)$ satisfies
the conditions of Theorem \ref{thmRG} and $G = G(x), x \in \R$  is  a Borel 
 function such that 
\begin{equation} \label{GMconv}
\lim_{M \to \infty} \E G\Big(\frac{X_{M}(\boldsymbol{0}) - \E X_M(\boldsymbol{0})}{M^{1/2}} \Big)^2 =  \E G(Z_\mu)^2  < \infty.
\end{equation}
Let $M  = \lambda^{\gamma} $ for some $\gamma > 0$.
Then
for any $\phi \in  \Phi $  as $\lambda \to \infty $
\begin{eqnarray}\label{limYinfty}
\lambda^{(\gamma/2) -H(\gamma)}(Y_{\lambda,M}(\phi) -
\E Y_{\lambda, M}(\phi)) &\limd&h_{G, \mu}(1)
\begin{cases}
B_\alpha(\phi),  &  \gamma > d(\alpha-1),  \\
L_{\alpha}(\phi), &   \gamma < d(\alpha-1),  \\
J_\alpha(\phi),  & \gamma = d(\alpha-1),
\end{cases}
\end{eqnarray}
where $H(\gamma), B_\alpha(\phi), L_\alpha(\phi), J_\alpha(\phi) $ are the same as in \eqref{limXinfty}, and
\begin{equation}\label{h1}
h_{G, \mu}(1) =   \mu^{-1} \E G(Z_\mu) Z_\mu
\end{equation}
is the first coefficient in the Hermite expansion  \eqref{Ghermite}  of $G$.

\end{theorem}

\noi {\bf Proof.} We have $Y_{M}(\mbt) = G_M (X_{M}(\mbt)) $ where
$G_M (x) := G( (x - \mu M)/M^{1/2}), x \in \N $ and $X_{M}(\mbt) =: N_M$ has a Poisson distribution  with mean $\mu M$.
Following \eqref{Gexp} consider the expansion
\begin{equation}\label{Gtauexp}
G_M(x) = \sum_{k=0}^\infty \frac{c_{G,M}(k)}{k!}  P_k(x; \mu M), \qquad x \in \N
\end{equation}
where
\begin{equation}\label{CGtau}
c_{G,M}(k) := (\mu M)^{-k} \E [ G_M (N_M ) P_k(N_M ; \mu M)],  \quad k \in \N
\end{equation}
Particularly,
\begin{eqnarray}\label{CG1conv}
M^{1/2} c_{G,M}(1)
&=&\frac{ M^{1/2}}{\mu M} \E\Big[ G\big(\frac{N_M - M \mu}{M^{1/2}}\big)(N_M - M \mu) \Big] \  \to  \    h_{G, \mu}(1)
\end{eqnarray}
as $M \to  \infty $, the limit as in  \eqref{h1}. The convergence in \eqref{CG1conv} under condition \eqref{GMconv}
can be verified with the help of Pratt's lemma \cite{pra1960}, as follows. Let $g (x) :=  G(x) x, x \in \R, \,
\xi_M := \frac{N_M - M \mu}{M^{1/2}}  $,
then $M^{1/2} c_{G,M}(1) =  \mu^{-1} \E g (\xi_M) \xi_M  $ and  $\xi_M \limd Z_\mu \, (M \to \infty)$. Since $G$ is a.e. continuous and
$Z_\mu$ has a continuous distribution, this implies
$G (\xi_M) \limd G(Z_\mu), \, g (\xi_M) $  $\limd g(Z_\mu) $.  By the Skorohod representation theorem,
$\xi_M \to_p Z_\mu,  G (\xi_M) \to_p G(Z_\mu), g (\xi_M) \to_p g(Z_\mu) $ in probability. Moreover, $ |g (\xi_M)| \le
(1/2)(G(\xi_M)^2 + \xi_M^2) =: \bar g_M $ where $\bar g_M \to_p (1/2)(G(Z_\mu)^2  $  $+ Z^2_\mu)  =: \bar g $ and
$\E \bar g_M \to \E \bar g  $ according to \eqref{GMconv}.  This and \cite{pra1960} prove \eqref{CG1conv}.

In view of \eqref{Gtauexp}, we have the representation
\begin{eqnarray}\label{Ydec}
Y_{\lambda,M}(\phi)- \E Y_{\lambda, M}(\phi)
&=&c_{G,M}(1) (X_{\lambda, M}(\phi) - \E X_{\lambda, M}(\phi)) +   Y^*_{\lambda,M}(\phi),  \qquad \text{where}\\
Y^*_{\lambda, M}(\phi)&:=&\int_{\R^d}  Y^*_{M}(\mbt) \phi(\mbt/\lambda) \d \mbt, \qquad
Y^*_{M}(\mbt):= \sum_{k=2}^\infty  \frac{c_{G,M}(k)}{k!}  P_k\big(X_M (\mbt); \mu M\big). \nn
\end{eqnarray}
The convergence in  \eqref{limXinfty} follows from \eqref{limYinfty} and \eqref{CG1conv}, once we show
that $Y^*_{\lambda, M}(\phi)$ in
\eqref{Ydec} is negligible, or
\begin{eqnarray}\label{Ystar}
\E Y^*_{\lambda, M}(\phi)^2&=&o(\lambda^{2H(\gamma)-\gamma}), \qquad \lambda \to \infty
\end{eqnarray}
for $M, H(\gamma)$ as in  Theorem  \ref{thmRGnon} and any fixed $\gamma  >0$.  Applying Corollary \ref{corM} with $k^* = 2,
\mu_3 = {\rm Cov}(X_M (\mbt), X_M (\boldsymbol{0})) = M  r_X (\mbt)  $
and the bound $c^2_{G,M}(k) \le (\mu M)^{-k} k! \E G_M ( X_M (\boldsymbol{0}))^2 $, see  \eqref{CGbdd},
we get
\begin{eqnarray}
|\E  Y^*_{M}(\mbt)  Y^*_{M}(\boldsymbol{0})|
&\le&\sum_{k=2}^\infty \frac{c^2_{G,M}}{(k!)^2} \big|\E [P_k(X_M (\mbt); \mu M) P_k(X_M (\boldsymbol{0}); \mu M)]  \big| \nn \\
&\le&\E G_M ( X_M (\boldsymbol{0}))^2  \sum_{k=2}^\infty \frac{ ( r_X (\mbt)M)^k }{ (\mu M)^k} \nn \\
&=&\E G_M ( X_M (\boldsymbol{0}))^2  \sum_{k=2}^\infty \big(\frac{r_X (\mbt) }{r_X(\boldsymbol{0})} \big)^k \ \le \ C (1 \wedge r^2_X(\mbt)).
\label{Ystar1}
\end{eqnarray}
Applying Proposition \ref{propcov} \eqref{c3}, relation \eqref{Ystar} follows since
$\max\{ d, 2d (2-\alpha)\} < \min \{ (3-\alpha) d, \frac{2(\gamma + d)}{\alpha} - \gamma \} $ holds for any
$\gamma >0, \alpha \in (1,2)$. \hfill $\Box$

\begin{rem} \label{remG} {\rm Condition \eqref{GMconv} on $G$ involving convergence of the second moments  only is rather weak. Using
the notation in \eqref{Ylambda}, it writes as
$\lim_{M \to \infty} \E Y_M(0)^2 =  \E G(Z_\mu)^2   < \infty. $
 \eqref{GMconv} can be replaced
by a boundedness condition:
\begin{eqnarray} \label{Gexp}
|G(x)| \le C_1 \e^{C_2 |x|}, \quad x \in \R \quad (\exists \, C_1, C_2 >0).
\end{eqnarray}
Indeed, verification of \eqref{GMconv} for $G(x) = C_1 \e^{C_2 |x|} $ is easy, implying   \eqref{GMconv} for $G$ in \eqref{Gexp} by
Pratt's lemma.
}
\end{rem}

\smallskip

Obviously, Theorem \ref{thmRGnon} does not hold when $\gamma = 0$ or $M = 1$ is fixed.
This case is treated in the following proposition.

\begin{proposition} \label{propGlim0} Let $Y(\mbt):= G\big(X(\mbt)\big), \mbt \in \R^d, $ \  where $X(\mbt) $ is as in \eqref{RG1} and
$G(x), x \in \N $ satisfies $\E G(X(\boldsymbol{0})^2 < \infty $,
$Y_{\lambda}(\phi):= \int_{\R^d}  Y(\mbt) \phi(\mbt/\lambda) \d \mbt,  \phi \in \Phi. $  Then
\begin{eqnarray} \label{limY0}
\lambda^{-d/\alpha} (Y_\lambda (\phi) - \E Y_\lambda (\phi))
&\limd&c_G(1;\mu) L_{\alpha}(\phi),
\end{eqnarray}
where $c_G(1;\mu) = \E G(X(\boldsymbol{0}))( X(\boldsymbol{0}) -  \E X(\boldsymbol{0})) $ and
$L_\alpha (\phi)$ is the same $\alpha$-stable RF as in \eqref{limX0}.

\end{proposition}

\noi {\bf Proof.} Similarly as in the proof of Theorem \ref{thmRGnon}, write
$Y_{\lambda}(\phi)- \E[Y_{\lambda}(\phi)]
= c_{G}(1;\mu) (X_{\lambda}(\phi) - \E [X_{\lambda}(\phi)]) + $  $  Y^*_{\lambda}(\phi), $ where
$Y^*_{\lambda}(\phi):= \int_{\R^d}  Y^*(\mbt) \phi(\mbt/\lambda) \d \mbt $   and
$Y^*(\mbt):= \sum_{k=2}^\infty  \frac{c_{G}(k;\mu)}{k!}  P_k\big(X(\mbt); \mu\big), $  $ \mu = \E [X(\boldsymbol{0})] $ satisfies
$|\E  Y^*(\mbt)  Y^*(\boldsymbol{0})| \le  C (1 \wedge r^2_X(\mbt)) $ as in \eqref{Ystar1}. Then, \eqref{limY0}
follows in view of  \eqref{limX0}  and \eqref{c3}. \hfill $\Box$

\begin{rem} \label{remran} {\rm Theorem \ref{thmRGnon} and Proposition \ref{propGlim0} yield
trivial limits if the respective coefficients  $h_{G, \mu}(1), $  $ c_G(1)$  vanish. The question of the limit distribution
of $Y_{\lambda,M}(\phi), Y_{\lambda}(\phi)$ in such case is open.
}
\end{rem}

\begin{example} \label{ex1}  (Scaling limit of the Boolean model. ) The Boolean model
$\hat X(\mbt)
= X(\mbt) \wedge 1 $
corresponds to $ Y(\mbt) = G( X(\mbt)) $ with $G(x) = x \wedge 1, x \in \N$. We have
$c_G(0;\mu) = 1 - \e^{-\mu},  c_G(k;\mu) = (-1)^{k+1} \e^{-\mu} \, (k \ge 1) $ and the convergence
in \eqref{limY0} holds with $c_G(1;\mu) = \e^{-\mu}$. Let $\phi (\mbx) = \I (x \in A), $  where
$A \subset \R^d $ is a Borel set and  $\hat X_\lambda (A)$  be as in \eqref{hatXA}.
From Proposition \ref{propGlim0}  and \eqref{limX0} it follows

\begin{corollary} \label{corBool} Let $A \subset \R^d $ be a bounded Borel set and $\hat X_\lambda(A) $ as in \eqref{hatXA}.
Then
\begin{equation}
\lambda^{-d/\alpha}(\hat X_\lambda (A) - \E \hat X_\lambda (A)) \limd \e^{-\mu}  L_\alpha (A), \quad \lambda \to \infty
\end{equation}
where $L_\alpha (A)$ is $\alpha$-stable r.v. with characteristic function $\E \e^{\i \theta L_\alpha (A)}
= \exp\{ - \sigma_\alpha |\theta|^\alpha \Leb_d (A)(1 - \i \, {\rm sgn}(\theta)\, \tan (\pi \alpha/2))\}, \theta \in \R. $
\end{corollary}

\end{example}

\begin{example}\label{ex2}  (Scaling limits of the Exponential RG model. )
We define the {\it Exponential RG model} as ${\cal E}(\mbt) := \e^{a X(\mbt)}, \mbt \in \R^d $ where $X(\mbt) $ is the RG model in \eqref{RG1} and
$ a \in \R $ a real parameter. We also consider the exponential function  of the aggregated RG model
\begin{eqnarray}\label{YYlambda}
{\cal E}_M (\mbt)&:=& \e^{a (X_M (\mbt) - \E X_M (\mbt))/M^{1/2}},  \qquad
{\cal E}_{\lambda, M}(\phi) := \int_{\R^d}  \phi(\mbt/\lambda) {\cal E}_M (\mbt) \d \mbt,
\end{eqnarray}
where $ X_M,  \phi$ are as in \eqref{Xlambda1}. As noted in the Introduction,  the interest in the scaling limits of
\eqref{YYlambda} is motivated by application to large-time asymptotics of statistical solution of Burgers' equation discussed in the last section.  Obviously, \eqref{YYlambda} is a particular case of
 \eqref{Ylambda} corresponding to $ G(x) = \e^{a x}$. Note $D^k_+ G(x) = (\e^a - 1)^k \e^{a x} $ and
$c_G (k) = (\e^a -1)^k \e^{(\e^a -1) \mu}, k \in \N$. We also  find that
\begin{eqnarray}
M^{1/2} c_{G,M}(1)
&=&\exp  \{ (\e^{a/M^{1/2}} - 1  - (a/M^{1/2})) \mu M \}M^{1/2} (\e^{a/M^{1/2}} - 1) \nn \\
&\to&a \e^{a^2 \mu /2}
\ = \ \mu^{-1} \E [\e^{a Z_\mu } Z_\mu] = h_{G,\mu}(1) \label{h1exp}
\end{eqnarray}
as $M \to \infty $, see \eqref{h1}. It is easy to see  $ G(x) = \e^{a x}$ that satisfies the conditions of
Theorem \eqref{thmRGnon} and the convergences in \eqref{limYinfty} hold for \eqref{YYlambda} with $h_{G,\mu}(1)$ in
\eqref{h1exp}. The relevant bound in \eqref{Ystar1} for the above $G(x)$ can be directly obtained from the equality
\begin{eqnarray*}
{\rm Cov}\big({\cal E}^*_M (\boldsymbol{0}), {\cal E}^*_M (\mbt) \big)
&=&(\E {\cal E}_M (\boldsymbol{0}))^2
\big\{\e^{(\e^{a/M^{1/2}} -1)^2 M r_{X}(\mbt)} -1
-(\e^{a/M^{1/2}} -1)^2 M r_{X}(\mbt)\big\},
\end{eqnarray*}
where ${\cal E}^*_M (\mbt) = ({\cal E}_M (\mbt) -  \E {\cal E}_M (\mbt))   - c_{G,M}(1)
(X_M (\mbt) -  \E X_M (\mbt)) $ as in \eqref{Ystar1}.
\end{example}

\begin{rem} \label{remsmall} {\rm \cite{bier2010} discusses {\it small-scale} scaling limits of RG model in \eqref{Xlambda1} as $\lambda \to 0 $ and
$M \to 0$ together with $\lambda$, under a similar condition (c.f. \eqref{falpha}) on the behavior of $f(r)$
as $r \to 0$ with $\alpha \in (0,1)$. Extending these results to nonlinear functions in \eqref{Ylambda}  is open.  Anisotropic
small-scale limits (without aggregation) for L\'evy driven RFs on $\R^2 $ were studied in \cite{pils2022}.

}
\end{rem}

\begin{rem} \label{remshot} {\rm  A class of RF which lie between Gaussian and RG RFs and are quite popular
in applied sciences are  {\it shot-noise}  RFs having a representation
\begin{equation}\label{SN}
X(\mbt) := \sum_{j=1}^\infty   W_j (\mbt - \mbu_j), \qquad \mbt \in  \R^d
\end{equation}
w.r.t. the same Poisson point process $\{\mbu_j \} $ as in \eqref{RG}, where
$W_j = \{ W_j(\mbt); \mbt \in \R^d \} $ are i.i.d. copies of (generic) {\it pulse  RF} $W = \{ W(\mbt); \mbt \in \R^d \} $, all independent of
$\{\mbu_j \} $.  We see that \eqref{SN} encompasses \eqref{RG} which  correspond to
$ W(\mbt)  = \I( \mbt \in R^{1/d} \Xi^0)$. Assuming that the trajectories of $W$ belong to $L^1 (\R^d)$ a.s.,
\eqref{SN} can be written as the Poisson stochastic integral
\begin{equation}\label{SN1}
X(\mbt) = \int_{\R^d \times L^1(\R^d)}   w(\mbt - \mbu) {\cal N}(\d \mbu, \d w), \qquad \mbt \in  \R^d
\end{equation}
w.r.t. to Poisson random measure  ${\cal  N}(\d u, \d w)$ with mean $ \E {\cal N}(\d \mbu, \d w) =  \d \mbu \P (W (\cdot) \in \d w), (\mbu, w)
\in   \R^d \times L^1(\R^d)$. The integral in \eqref{SN1} is a well-defined and stationary RF with finite variance provided
$\int_{\R^d} (\E |W(\mbt)| + \E |W(\mbt)|^2)  \d \mbt < \infty $ holds, in which case $\E X(\mbt) = \int_{\R^d} \E W(\mbt) \d \mbt$
and the covariance
${\rm Cov}(X (\boldsymbol{0}), X(\mbt))
= \int_{\R^d} \E W(\mbu) W(\mbt + \mbu) \d \mbu, \ \mbt \in \R^d $
may exhibit LRD property under suitable assumption on the pulse RF $W$. A rather general form of pulse allowing
for LRD (see \cite{alb1994, gir1992, surw1994})  is given by
\begin{eqnarray} \label{SN4}
W(\mbt) = \eta \, a(R^{-1/d} \mbt), \qquad \mbt \in \R^d,
\end{eqnarray}
where $a (\cdot) \in L^2 (\R^d)$ is a deterministic function whereas $\eta \in \R, R >0$ are r.v.s satisfying certain moment conditions.
The  corresponding covariance function writes as
 \begin{eqnarray}\label{SN5}
{\rm Cov}(X (\boldsymbol{0}), X(\mbt))
&=&\E [\eta^2 R \, (a \star a)\big(R^{-1/d} \mbt\big)\big],  \qquad \mbt \in \R^d,
\end{eqnarray}
where $\star $ denotes convolution. \eqref{SN4} may be interpreted as a `typical'   `$a$-pulse' with random `amplitude' $\eta $ and random `frequency' $R^{-1/d} $. \cite[Lemma 2]{gir1992} provides general conditions on $\eta, R, a\star a $ implying covariance LRD and regular decay of \eqref{SN5} as
$\|\mbt\|\to \infty$ for $d=1$ which can be generalized to any $d \ge 1 $.  Gaussian scaling limits of shot-noise RFs in \eqref{SN}, \eqref{SN4}
were studied in several papers (see, e.g., \cite{gir1992} and the references therein). The recent work \cite{leip2022} discussed scaling limits
and a trichotomy similar to  \eqref{limXy} for one-dimensional ($d=1$)  shot-noise in \eqref{SN}  with rescaled intensity
$M \d \mbu \P (W (\cdot) \in \d w), M = \lambda^\gamma $ of the Poisson random measure, for a general class of pulse $W$. An interesting
open problem is to  extend the results on nonlinear functionals in Theorem   \ref{thmRGnon} to shot-noise RF  in \eqref{SN1}-\eqref{SN4}.

}
\end{rem}

\begin{rem} \label{remcox} {\rm A {\it Cox} (or doubly stochastic Poisson) point process $ U = \{ \mbu_j\} $  is a Poisson point process on $\R^d $
with {\it random} intensity $\zeta (\mbu) \d \mbu $, where $ \{ \zeta (\mbu); \mbu \in \R^d \} =: \zeta $
is a nonnegative RF, meaning that  the conditional distribution of $U$ given $\zeta $ is Poisson with intensity  $\zeta (\mbu) \d \mbu $ \cite{sto1989}.
Scaling limits of shot-noise RFs driven by Cox process were studied in
\cite{funa1995} and other works. LRD property in such RFs may be due to random intensity $\zeta $. Extending some results of this sec.
to RG models driven by Cox process seems feasible.

}
\end{rem}

\section{Application: scaling limits of solutions of Burgers' equation with initial
RG data}

Burgers' equation with (random) potential initial data is written as
\begin{eqnarray}\label{burg}
\partial \vec v(t, \mbx)/\partial t + (\vec v(t,\mbx), \nabla) \vec v(t,\mbx)
&=&{\small \frac{1}{2}} \kappa \Delta \vec v(t, \mbx), \\
\vec v(0, \mbx)&=&- \nabla \xi (\mbx), \nn
\end{eqnarray}
where $\vec v(t,\mbx) = (v_1(t, \mbx), \cdots, v_d (t, \mbx)), (t, \mbx) \in \R_+ \times \R^d $ is a $\R^d$-valued function (velocity field) and
$\xi =  \{\xi (\mbx); \mbx \in \R^d \} $ is a scalar (potential) RF;
$(\vec v(t,\mbx), \nabla) := \sum_{i=1}^d v_i(t, \mbx) \partial/\partial x_i$.
The parameter $\kappa >0 $ is usually called the viscosity parameter.
Burgers' equation is one of the important equations of mathematical physics. The solution $ \vec v(t,\mbx)$ with random initial data is a (vector-valued) RF whose behavior as $t \to \infty $ and/or $\kappa \to 0 $ presents considerable physical and mathematical  interest and has been extensively studied
in the literature. For a probabilistic approach, we refer to \cite{alb1994} and the review paper \cite{surw2003}. When $\kappa >0$ is fixed,  the natural
parabolic scaling leads to  the RF $\vec v_\lambda(t, \mbx) := \vec v(\lambda^2  t,  \lambda \mbx) $ and the problem concerns
the limit distribution of RF $\vec v_\lambda(t, \mbx); (t, \mbx) \in \R_+ \times \R^d $ as $\lambda \to \infty $.

The study of Burgers' equation is facilitated by the Hopf-Cole substitution
\begin{equation}\label{vecv}
\vec v(t, \mbx) = - \kappa \nabla \log u(t, \mbx) = - \frac{\kappa \nabla u(t, \mbx)}{u(t, \mbx)}
\end{equation}
with a scalar-valued $u(t, \mbx)$ satisfying the heat equation $\partial u(t, \mbx)/\partial t
= {\small \frac{1}{2}} \kappa \Delta u(t, \mbx) $ with the (exponential) initial condition $u(0+, \mbx) = \e^{\xi (\mbx)/\kappa}, \mbx \in \R^d.   $
Thus, \eqref{vecv} has an explicit representation through the heat kernel $g(t, \mbx, \mby) := (2\pi \kappa t)^{-d/2}
\exp \{ - \|\mbx - \mby\|^2/2 \kappa t\} $ as the ratio
\begin{equation}\label{vecv1}
\vec v(t, \mbx) = -\frac{ \kappa \int_{\R^d}   \nabla g(t, \mbx, \mby) \e^{ \xi(\mby) /\kappa } \d \mby }
{\int_{\R^d}  g(t, \mbx, \mby) \e^{ \xi(\mby) /\kappa } \d \mby }.
\end{equation}
Using the fact that $\int_{\R^\d}  \nabla g(t, \mbx, \mby) \d \mby = 0$, one can replace
$\e^{ \xi(\mby) /\kappa }$ in the numerator of \eqref{vecv1} by $\e^{ \xi(\mby) /\kappa } - \E \e^{ \xi(\mby) /\kappa } $,
provided the last expectation is finite and does not depend on $\mby$. As a consequence, the  rescaled velocity RF writes as
\begin{equation}\label{vecv2}
\vec v_\lambda (t, \mbx) = -\frac{\kappa \int_{\R^d} \phi_{t, \mbx} (\mby/\lambda) (G(\xi (\mby))  - \E  G(\xi(\mby))) \d \mby}
{ \lambda \int_{\R^d} \psi_{t, \mbx} (\mby/\lambda) G(\xi(\mby)) \d \mby},
\end{equation}
where $G(x) = \e^{x/\kappa}$ and
the integrals in the  numerator and denominator resemble \eqref{Xlambda0} with $\phi (\mby) = \phi_{t, \mbx}(\mby)
:= \nabla g(t, \mbx, \mby) $ and
 $\phi (\mby) = \psi_{t,\mbx}(\mby)  := g(t, \mbx, \mby)$, respectively. Clearly, for any fixed $(t,  \mbx)  \in \R_+ \times \R^d $ the above $\phi$'s belong to
$ \Phi$ in  \eqref{Phi}.  

As mentioned in the Introduction, our  aim is the limit distribution of \eqref{vecv2} for the initial potential RF
\begin{equation}\label{xitau}
\xi (\mby) := M^{-1/2}(X_M (\mby) - \E X_M (\mby)), \quad \mby \in \R^d,
\end{equation}
where $X_M $ is the (aggregated) RG model as in Theorem \ref{thmRGnon} and elsewhere in this paper, with intensity $M =  \lambda^\gamma $ increasing
with $\lambda $ when $\gamma >0$. Similarly as in \cite{alb1994} and some other works, we ignore the meaning of the initial condition
in \eqref{burg}  for $\xi$ in \eqref{xitau} and consider scaling properties of the Hopf-Cole solution in  \eqref{vecv1} alone. Note
$\E \e^{\xi(\mby)/\kappa} \to \E \e^{Z_\mu/\kappa} = \e^{\mu/2\kappa^2} \,  (M\to  \infty),
\E [\e^{Z_\mu/\kappa} Z_\mu]/\E \e^{Z_\mu/\kappa}= 1/\kappa $ and
$\lambda^{-d}  \int_{\R^d}  g(t, \mbx, \mby/\lambda) \e^{ \xi(\mby) /\kappa } \d \mby \to \E \e^{Z_\mu/\kappa} $  in probability by the law of large numbers.
The application of  Theorem \ref{thmRGnon} with $G(x) = \e^{x/\kappa}, x \in \R$ and \eqref{h1exp} (Example 2) yields
the following result.

\begin{corollary} Let $\vec v_\lambda (t, \mbx)$ be as in \eqref{vecv2}, \eqref{xitau},  with $X_M, M  $ satisfying the conditions
of Theorem \ref{thmRGnon}. Then, as $\lambda \to \infty $,  for any $\gamma >0$
\begin{eqnarray} \label{burlim}
\lambda^{1+ d + \frac{\gamma}{2} - H(\gamma)} \vec v_\lambda(t, \mbx)
&\limfdd&
\begin{cases}
B_\alpha(\nabla g(t,\mbx, \cdot)),  &  \gamma > d(\alpha-1),  \\
L_{\alpha}(\nabla g(t,\mbx, \cdot)), &   \gamma < d(\alpha-1),   \\
J_\alpha(\nabla g(t,\mbx, \cdot)),  & \gamma = d(\alpha-1),
\end{cases}
\end{eqnarray}
where $H(\gamma)$ and the limit RFs are the same as in \eqref{limXinfty}

\end{corollary}

Let us remark that convergence to $\alpha$-stable limit in \eqref{burlim} holds also for $\gamma = 0$ or
$\xi(\mby)$ in \eqref{xitau} replaced by $\xi(\mby) = X(\mby) $ in \eqref{RG1}.  From
Proposition \ref{propGlim0} we conclude the following result.

\begin{corollary} Let $\vec v_\lambda (t, \mbx)$ be as in \eqref{vecv2}  with $\xi(\mby) = X(\mby) $ given in \eqref{RG1}.
Then, as $\lambda \to \infty $
\begin{eqnarray} \label{burlim0}
\lambda^{1+ d  -  \frac{d}{\alpha}} \vec v_\lambda(t, \mbx)
&\limfdd&\kappa (\e^{1/\kappa} -1)
L_{\alpha}(\nabla g(t,\mbx, \cdot)),
\end{eqnarray}
where $L_\alpha $ is defined in  \eqref{Lalpha}-\eqref{JJ}.
\end{corollary}

We remark that in dimension $d=1$, a similar result to \eqref{burlim0} was proved in \cite[Thm.1.1(iii)]{sur1996} for initial  (potential) data $\xi(y), y \in \R$
given by a piecewise-constant renewal-reward process with renewal distribution having $\alpha$-tail with $\alpha \in (1,2)$.

\section*{Acknowledgement}

The author is grateful to two anonymous referees for useful  and constructive comments.

\end{document}